\input amstex\documentstyle {amsppt}  
\pagewidth{12.5 cm}\pageheight{19 cm}\magnification\magstep1
\topmatter
\title Notes on character sheaves\endtitle
\author G. Lusztig\endauthor
\address Department of Mathematics, M.I.T., Cambridge, MA 02139\endaddress
\dedicatory{To Pierre Deligne on the occasion of his 65-th birthday}\enddedicatory 
\thanks Supported in part by the National Science Foundation\endthanks
\endtopmatter   
\document
\define\tit{\ti\tau}
\define\tcw{\ti{\cw}}
\define\sgn{\text{\rm sgn}}
\define\ubK{\un{\bar K}}
\define\bfZ{\bar{\fZ}}
\define\bft{\bar{\ft}}
\define\bfs{\bar{\fs}}
\define\ubbq{\un{\bbq}}
\define\codim{\text{\rm codim}}
\define\De{\Delta}
\define\Lo{\text{\rm Loc}}
\define\uLo{\un{\Lo}}

\define\uA{\un A}

\define\uK{\un K}

\define\uL{\un L}

\define\hD{\hat D}

\define\da{\dagger}

\define\po{\text{\rm pos}}

\define\si{\sim}
\define\wt{\widetilde}
\define\sqc{\sqcup}

\define\bK{\bar K}

\define\bY{\bar Y}

\define\bZ{\bar Z}

\define\op{\oplus}

\define\part{\partial}
\define\em{\emptyset}

\define\m{\mapsto}
\define\do{\dots}

\define\bsl{\backslash}

\define\sub{\subset}    

\define\T{\times}
\define\ti{\tilde}
\define\nl{\newline}
\redefine\i{^{-1}}

\define\un{\underline}
\define\ov{\overline}
\define\ot{\otimes}
\define\bbq{\bar{\QQ}_l}

\define\Hom{\text{\rm Hom}}

\define\tr{\text{\rm tr}}

\define\supp{\text{\rm supp}}

\define\a{\alpha}
\redefine\b{\beta}
\redefine\c{\chi}
\define\g{\gamma}
\redefine\d{\delta}

\define\et{\eta}
\define\io{\iota}
\redefine\o{\omega}
\define\p{\pi}
\define\ph{\phi}

\define\r{\rho}

\redefine\t{\tau}

\redefine\l{\lambda}
\define\z{\zeta}
\define\x{\xi}

\redefine\G{\Gamma}
\redefine\D{\Delta}

\define\Th{\Theta}
\redefine\L{\Lambda}
\define\Ph{\Phi}
\define\Ps{\Psi}

\define\dd{\bold d}

\define\kk{\bold k}

\define\qq{\bold q}

\redefine\AA{\bold A}

\define\FF{\bold F}

\define\II{\bold I}

\define\NN{\bold N}

\define\QQ{\bold Q}

\define\WW{\bold W}
\define\ZZ{\bold Z}

\define\ca{\Cal A}
\define\cb{\Cal B}

\define\cd{\Cal D}
\define\ce{\Cal E}

\define\ch{\Cal H}

\define\ck{\Cal K}
\define\cl{\Cal L}

\define\cs{\Cal S}
\define\ct{\Cal T}

\define\cw{\Cal W}
\define\cz{\Cal Z}

\define\fc{\frak c}
\define\fd{\frak d}
\define\fe{\frak e}

\define\fs{\frak s}
\define\ft{\frak t}

\define\fC{\frak C}
\define\fD{\frak D}

\define\fK{\frak K}

\define\fZ{\frak Z}

\define\tc{\ti c}

\define\tB{\ti B}

\define\tY{\ti Y}
\define\tZ{\ti Z}

\define\sh{\sharp}

\define\tce{\ti\ce}
\define\bul{\bullet}

\define\ucl{\un\cl}

\define\BV{BV}
\define\BFO{BFO}
\define\BBD{BBD}
\define\GI{Gi}
\define\KL{KL1}
\define\KLL{KL2}
\define\MV{MV}
\define\CS{L1}
\define\RA{L2}
\define\PAR{L3}
\define\CDG{L4}

\head Introduction\endhead
Let $\kk$ be an algebraically closed field and let $G$ be an affine algebraic group over $\kk$ which is reductive
(that is, its identity component $G^0$ is reductive). Let $\NN_\kk^*$ be the set of all integers $\ge1$ that are 
$\ne0$ in $\kk$. We fix a prime number $l\in\NN_\kk^*$. According to the theory of character sheaves (see 
\cite{\CS, I} for $G=G^0$ and \cite{\CDG, VI} in the general case) one can define a natural class of simple 
perverse $\bbq$-sheaves on $G$ (the "character sheaves" on $G$) whose properties mimic those of the irreducible 
characters of a reductive group over a finite field.  

This note has two parts. In \S1 we study the functor $\b:=\ft_*\fs^!\fc_!\fd^*$ (see below) introduced by 
Bezrukavnikov, Finkelberg and Ostrik \cite{\BFO} (who assumed that $G=G^0$ and that the characteristic of $\kk$ is
$0$); they prove the remarkable result that the complex obtained by the application of $\b$ to a character sheaf 
is a perverse sheaf. Here the following notation is used:

$D$ is a fixed connected component of $G$;

$\cb$ is the variety of Borel subgroups of $G^0$; for any $B\in\cb$, $U_B$ is the unipotent radical of $B$;

$\tZ=\{(B,B',g)\in\cb\T\cb\T D;gBg\i=B'\}$;

$Z=\{(B,B',gU_B);(B,B')\in\cb^2,g\in D;gBg\i=B'\}$;

$\fZ=\{(B_1,B,B',gU_B);(B_1,B,B')\in\cb^3,g\in D,gBg\i=B';B_1,B\text{ opposed}\}$;

$Z'=\{(B_1,B',U_{B'}gU_{B_1});(B_1,B')\in\cb^2,g\in D;
gB_1g\i,B'\text{ opposed}\}$;

$\fd:\tZ@>>>D$ is $(B,B',g)\m g$;

$\fc:\tZ@>>>Z$ is $(B,B',g)=(B,B',gU_B)$;

$\fs:\fZ@>>>Z$ is $(B_1,B,B',gU_B)\m(B,B',gU_B)$;

$\ft:\fZ@>>>Z'$ is $(B_1,B,B',gU_B)\m(B_1,B',U_{B'}gU_{B_1})$;

$\fd^*:\cd(D)@>>>\cd(\tZ)$, $\fc_!:\cd(\tZ)@>>>\cd(Z)$, $\fs^!:\cd(Z)@>>>\cd(\tZ)$, $\ft_*:\cd(\tZ)@>>>\cd(Z')$ 
are the corresponding Grothendieck functors.
\nl
(For an algebraic variety $X$ over $\kk$ we write $\cd(X)$ for the bounded derived category of $\bbq$-sheaves on 
$X$.) Actually in \cite{\BFO} (where $G=G^0$), the varieties $Z,Z'$ appear in a different (but equivalent) form as
$T\bsl(G/U_B\T G/U_B)$, $T\bsl(G/U_B\T G/U_{B'})$, where $B,B'$ are two opposed Borel subgroup and $T=B\cap B'$.

Note that (in 1987) I showed that the cohomology sheaves of $\fc_!\fd^*(A)$ for a character sheaf $A$ on $D$ (when
$D=G^0$) have a particularly simple behaviour and that this behaviour characterizes character sheaves (see 
\cite{\MV}, \cite{\GI}). On the other hand, the functor $\ft_*\fs^!$ is essentially an intertwining operator.

In the remainder of this paper we assume that $\kk$ is an algebraic closure of a finite field 
$\FF_q$.

In \S1 we restrict ourselves for simplicity to the case of unipotent character sheaves. (In some respects general
character sheaves behave like unipotent character sheaves on a possibly smaller group.) 
The main result in \S1 is an explicit computation of $\b(A)$ in a Grothendieck group which 
takes weights into account (under a mild restriction on the characteristic of $\kk$).
See 1.2(b) for a precise statement.

We now describe the content of \S2. We would like to understand how the tensor product of two irreducible 
representations $\r,\r'$ (over $\bbq$) of a reductive group $\G$ over $\FF_q$ decomposes into irreducibles. Take 
for example $\G=PGL_2(\FF_q)$. Let $\Th$ be the $(q^2-1)$ dimensional representation in which each irreducible 
representation of $\G$ (other than the unit representation) appears exactly once. If $\r,\r'$ are two irreducible
constituents of $\Th$ then $\r\ot\r'$ is equal (as a virtual representation) to $\Th$ plus or minus the sum of at
most three irreducible representations. (If $\r=\r'$ is the Steinberg representation then $\r\ot\r'$ is $\Th$ plus
the unit representation; if $\r$ and $\r'$ are two principal series representations then, most of the time,
$\r\ot\r'$ is $\Th$ plus a sum of two principal series representations; if $\r$ and $\r'$ are two discrete series
representations then, most of the time, $\r\ot\r'$ is $\Th$ minus a sum of two discrete series representations; if
$\r$ is a principal series representation and $\r'$ is a discrete series representation then $\r\ot\r'=\Th$.) We 
see that while the character of $\r,\r'$ can be described in terms of character sheaves on $PGL_2(\kk)$, the 
character of $\r\ot\r'$ cannot be described in terms of character sheaves (due to the presence of $\Th$). Note that the character of $\Th$ is the function with value $-1$ at regular unipotent elements, value $q^2-1$ at $1$ and value $0$ elsewhere. This is a linear combination of two class
functions on $\G$ which are characteristic functions of two simple perverse sheaves on 
$PGL_2(\kk)$ (one supported by the unipotent variety and one supported
by the unit element). If we enlarge the class of character sheaves by including these two
simple perverse sheaves the resulting class of simple perverse sheaves has the property
that the tensor product of two members in the class is a suitable combination of members of the
class, unlike the (unenlarged) class of character sheaves.
In \S2 we show how to enlarge (for a general $G$) the class of character sheaves to 
a larger class of simple perverse sheaves with a similar behaviour under tensor product
as in the case of $PGL_2(\kk)$.

{\it Notation.} We shall use extensively the notation and results of \cite{\BBD}. For $s\in\ZZ_{>0}$ let 
$\FF_{q^s}$ be the subfield of $\kk$ of cardinal $q^s$. Let $\ca=\ZZ[v,v\i]$ ($v$ an indeterminate). 

If $X$ is an algebraic variety over $\kk$ and $K\in\cd(X)$, $n\in\ZZ$, we write $K[[n]]$ instead of $K[n](n/2)$. 
Let $\fD_X:\cd(X)@>>>\cd(X)$ denote Verdier duality. For $K\in\cd(X)$ let $\ch^i(K)$ be the $i$-th cohomology 
sheaf of $K$ and let $\ch^i(K)_x$ be its stalk at $x\in X$; we write $H^iK$ instead of ${}^pH^iK$. If $X'$ is a 
closed subvariety of $X$, for any $K\in\cd(X')$ we set $K^X=j_!K\in\cd(X)$ where $j:X'@>>>X$ is the inclusion.

If $X$ has a given $\FF_q$-structure we write $\cd_m(X)$ for the corresponding mixed derived category of $\bbq$-sheaves. If 
$A\in\cd_m(X)$ is perverse and $j\in\ZZ$, we denote by $A_j$ the canonical subquotient of $A$ wich is pure of 
weight $j$. We write $X_s$ instead of $X(\FF_{q^s})$ and we denote by $\bbq^{X_s}$ the $\bbq$-vector space 
consisting of all functions $X_s@>>>\bbq$. If $K\in\cd_m(X)$ and $s\in\ZZ_{>0}$, we define a function 
$\c_{K,s}:X_s@>>>\bbq$ by
$$\c_{K,s}(\x)=\sum_{i\in\ZZ}(-1)^i\tr(F^s,\ch^i(K)_\x)$$
where $F$ is the Frobenius map relative to $\FF_q$.

Let $\WW$ be the set of $G^0$-orbits on $\cb\T\cb$ for the $G^0$-action given by conjugation on both factors. For
$B,B'\in\cb$ we write $\po(B,B')=w$ if the $G^0$-orbit of $(B,B')$ is $w$. We regard $\WW$ as a finite Coxeter 
group with length function $l:\WW@>>>\NN$ as in \cite{\CDG, 26.1}; let $\II=\{w\in\WW;l(w)=1\}$. Let $w_0$ be the 
longest element of $\WW$. Let $\le$ be the standard partial order on $\WW$. 

We shall assume that on $G$ we are given an $\FF_q$-structure with Frobenius map $F$ compatible with the group 
structure such that $F$ acts as identity on $G/G^0$ and on $\WW$.

For $g\in G$ let $Z_G(g_s)^0$ be the identity component of the centralizer in $G$ of the semisimple part of $g$.
Let $\cz^0_{G^0}$ be the identity component of the centre of $G^0$. For a connected component $D'$ of $G$ let 
${}^{D'}\cz^0_{G^0}$ be the set of all $g\in\cz^0_{G^0}$ such that $g$ commutes with some/any element of $D'$. If
$G_1$ is a subgroup of a group $G_2$ let $N_{G_2}(G_1)$ be the normalizer of $G_1$ in $G_2$.

{\it Acknowledgement.} I wish to thank the Institute for Advanced Study for its hospitality during April 2008 when
most of this work was done. I also whish to thank R. Bezrukavnikov for explaining to me some aspects of 
\cite{\BFO}.

\head 1. Study of the functor $\b$\endhead
\subhead 1.1\endsubhead
Let $X$ be an algebraic variety over $\kk$. Let $\ck(X)$ be the Grothendieck group of the category of perverse 
sheaves on $X$; it has $\ZZ$-basis given by the isomorphism classes of simple perverse sheaves on $X$. Let 
$\ck_\ca(X)=\ca\ot\ck(X)$. If $A,K$ are perverse sheaves on $X$ with $A$ simple we denote by $(A:K)$ the 
multiplicity of $A$ in a Jordan-H\"older series for $K$. Any perverse sheaf $K$ on $X$ gives rise to an element 
$\sum_A(A:K)A\in\ck(X)$ ($A$ runs over the isomorphism classes of simple perverse sheaves on $X$); this element is
denoted again by $K$. 

We define a symmetric bilinear inner product $(:):\ck_\ca(X)\T\ck_\ca(X)@>>>\ca$ by $(A:A')=1$ (resp. $(A:A')=0$)
if $A,A'$ are isomorphic (resp. nonisomorphic) simple perverse sheaves on $X$. If $A,K$ are perverse sheaves on 
$X$ with $A$ simple then the multiplicity $(A:K)$ and the inner product $(A:K)$ coincide. 

If $X$ has a given $\FF_q$-structure and $K\in\cd_m(X)$, we set 
$$gr(K)=\sum_{i,j\in\ZZ}(-1)^iv^jH^i(K)_j\in\ck_\ca(X).$$
Let $X,Y$ be algebraic varieties defined over $\FF_q$ and let $f:X@>>>Y$ be a morphism defined over $\FF_q$. We 
define a linear map $f^*_s:\bbq^{Y_s}@>>>\bbq^{X_s}$ by $f^*_s(\ph)(\x)=\ph(f(\x))$ for any $\ph\in\bbq^{Y_s}$, 
$\x\in X_s$. We define a linear map $f_{!s}:\bbq^{X_s}@>>>\bbq^{Y_s}$ by
$$f_{!s}(\ph')(\x')=\sum_{\x\in X_s;f(x)=\x'}\ph'(\x)$$
for any $\ph'\in\bbq^{X_s}$, $\x'\in Y_s$. Following Grothendieck we note that if $K\in\cd_m(X)$ then 
$f_!K\in\cd_m(Y)$ and we have $\c_{f_!K,s}=f_{!s}(\c_{K,s})$; if $K'\in\cd_m(Y)$ then $f^*K'\in\cd_m(X)$ and we 
have $\c_{f^*K',s}=f^*_s(\c_{K',s})$.

\subhead 1.2\endsubhead
For $w\in\WW$ let 

$\bZ^w=\{(B,B',xU_B)\in Z;\po(B,B')\le w\}$, 

$Z^w=\{(B,B',xU_B)\in Z;\po(B,B')=w\}$
\nl
(an open dense smooth irreducible subvariety of $\bZ^w$). Let 

$\bZ'{}^w=\{(B_1,B',U_{B'}xU_{B_1})\in Z';\po(B_1,B')\le w\}$, 

$Z'{}^w=\{(B_1,B',U_{B'}xU_{B_1})\in Z';\po(B_1,B')=w\}$ 
\nl
(an open dense smooth irreducible subvariety of $\bZ'{}^w$). 

Let $L^w=IC(\bZ^w,\bbq)^Z\in\cd(Z)$ where $\bbq$ is regarded as a local system on $Z^w$. Let $\bbq^w\in\cd(Z)$ be the extension by $0$ of the local system $\bbq$ on $Z_w$. Let 
$L'{}^w=IC(\bZ'{}^w,\bbq)^{Z'}\in\cd(Z')$ where $\bbq$ is regarded as a 
local system on $Z'{}^w$. We set

 $\De=\dim G^0=\dim D$, $M_w=\dim Z^w=\dim Z'{}^w=\De-l(w_0w)$
\nl
 so that 
$L^w[[M_w]]$ is a simple perverse sheaf on $Z$ and $L'{}^w[[M_w]]$ is a simple perverse sheaf on $Z'$. Let 
$$K^w_D=\fd_!\fc^*(\bbq^w)\in\cd(D),\bK^w_D=\fd_!\fc^*(L^w)\in\cd(D).$$
A simple perverse sheaf $A$ on $D$ is said to be a {\it unipotent character sheaf} if $(A:H^i(K^w_D))\ne0$ for
some $w\in\WW,i\in\ZZ$ or equivalently if $(A:H^i(\bK^w_D))\ne0$ for some $w\in\WW,i\in\ZZ$. Let $\hD^{un}$ be the
class of unipotent character sheaves on $D$. Let $\ck^{un}_\ca(D)$ be the $\ca$-submodule of $\ck_\ca(D)$ spanned
by the unipotent character sheaves on $D$. Let $\Xi$ be a set of representatives for the isomorphism classes of 
objects in $\hD^{un}$; note that $\Xi$ is a finite set. For any $A\in\Xi$ we set $d_A=\dim\supp(A)$, 
$d'_A=\codim_D\supp(A)$. If $A\in\Xi$ then $\fD_D(A)\in\hD^{un}$ and we denote by $A^*$ the object of $\Xi$ which
is isomorphic to $\fD_D(A)$. Under the $\ca$-linear involution $\dd:\ck^{un}_\ca(D)@>>>\ck^{un}_\ca(D)$ (the 
"duality" in \cite{\CDG, IX, \S42}), for any $A\in\Xi$, we have $\dd(A)=(-1)^{d'_A}A^\bul$ where $A^\bul$ is a 
well defined object of $\Xi$; moreover we have $d_A=d_{A^\bul}$.

We shall assume that the given $\FF_q$ structure on $G$ is such that and each $A\in\Xi$ satisfies $F^*A\cong A$.
(Such an $\FF_q$-structure exists since $\Xi$ is finite.) For each $A\in\Xi$ we can find (and we fix) an object 
$\uA\in\cd_m(D)$ which gives rise to $A$ and such that for any $g$ in an open dense subset of $\supp(A)$ and any 
$s\in\ZZ_{>0}$ such that $F^s(g)=g$, the eigenvalues of $F^s$ on the stalk $\ch^{-d_A}(A)_g$ are roots of $1$ 
times $q^{sd'_A/2}$. Note that $\uA$ is pure of weight $\De$. Now each of the varieties $\tZ,Z,\fZ,Z'$ has a 
natural $\FF_q$-structure induced by that of $G$ and the maps $\fd,\fc,\fs,\ft$ are defined over $\FF_q$. Hence we
have naturally $\b(\uA)\in\cd_m(Z')$ ($\b$ as in \S0) and $H^i(\b(uA))_j$ is well defined for any $i,j\in\ZZ$.

In the remainder of \S1 we shall make the following assumption:

(a) {\it either the characteristic of $\kk$ is a good prime for $G^0$ or $G^0$ is of classical type.}

\proclaim{Proposition} For any $A\in\Xi$ and any $j\in\ZZ$ we have
$$\sum_i(-1)^iH^i(\b(\uA))_j=\sum_{x\in\WW}(A^\bul:H^{j+l(w_0x)}(\bK^{w_0x}_D))L'{}^x[[M_x]]\tag b$$
in $\ck(Z')$. 
\endproclaim
The proof is given in 1.11.

Note that the right hand of (b) is the class in $\ck(Z')$ of a perverse sheaf on $Z'$. This suggests that 
$\b(\uA)$ is a perverse sheaf; by \cite{\BFO}, this is actually the case if the characteristic of $\kk$ is large 
enough.

\subhead 1.3\endsubhead
Let $H$ be the Iwahori-Hecke algebra attached to $\WW$ that is, the free $\ca$-module with basis $\{T_w;w\in\WW\}$
and with $\ca$-algebra structure given by $T_wT_{w'}=T_{ww'}$ if $l(ww')=l(w)+l(w')$, $(T_s+1)(T_s-v^2)=0$ if 
$s\in\II$. Note that $T_w$ is invertible in $H$ for any $w\in\WW$. Define an $\ca$-linear map $h\m{}^th$, $H@>>>H$
by ${}^tT_w=T_{w\i}$ for all $w$ (an algebra antiautomorphism). Define an $\ca$-linear map $h\m h{}^\da$, $H@>>>H$
by $T_w^\da=(-v^2)^{l(w)}T_{w\i}\i$ for all $w$ (an algebra involution commuting with $h\m{}^th$). We have a ring
involution $h\m\bar h$, $H@>>>H$ such that $v^jT_w\m v^{-j}T_{w\i}\i$ for all $w\in\WW,j\in\ZZ$. For $w\in\WW$ let
$$c_w=v^{-l(w)}\sum_{y\in\WW}P_{y,w}(v^2)T_y\in H,$$
where $P_{y,w}$ are the polynomials in the indeterminate $\qq$ defined in \cite{\KL}. We have $P_{y,w}=0$ unless 
$y\le w$. Moreover, $\ov{c_w}=c_w$. According to \cite{\KL}, the matrix $(P_{y,w})$ indexed by $\WW\T\WW$ has an 
inverse $(Q_{y,w})$ where 

(a) $Q_{y,w}=(-1)^{l(w)-l(y)}P_{w_0w,w_0y}$.

\subhead 1.4\endsubhead
For $x,y\in\WW$ we define $a_{x,y}\in\ca$, $b_{x,y}\in\ca$ by
$$T_{w_0}\i c_x=\sum_ya_{x,y}c_y,  (-v^2)^{-l(w_0)}T_{w_0}c_x=\sum_yb_{x,y}c_y.\tag a$$
In this subsection we prove for any $x,z\in\WW$ that:
$$a_{w_0x,w_0z}=a_{xw_0,zw_0}=(-1)^{l(x)-l(z)}b_{z,x}.\tag b$$
Let $\io:H@>>>H$ be the algebra involution of $H$ such that $\io(T_w)=T_{w_0ww_0}$ for all $w\in\WW$. Note that 
$\io(c_w)=c_{w_0ww_0}$ for all $w\in\WW$. Applying $\io$ to the first equation in (a) we obtain 
$T_{w_0}\i c_{w_0xw_0}=\sum_ya_{x,y}c_{w_0yw_0}$. On the other hand we have 
$T_{w_0}\i c_{w_0xw_0}=\sum_ya_{w_0xw_0,w_0yw_0}c_{w_0yw_0}$. It follows that $a_{w_0xw_0,w_0yw_0}=a_{x,y}$ for 
all $x,y$. This proves the first equality in (b).

To prepare for the proof of the second equality in (b), we set $H^*=\Hom_\ca(H,\ca)$ and we define a basis 
$(\tc_w)_{w\in\WW}$ of $H^*$ by $\tc_x(c_y)=(-1)^{l(x)}\d_{x,y}$. Define an $H$-module structure on $H$ by the 
left multiplication. For $h\in H,\ph\in H^*$ define $h\ph\in H^*$ by $(h\ph)(h_1)=\ph(hh_1)$ for $h_1\in H$; we 
define $h*\ph\in H^*$ by $h*\ph=({}^th^\da)\ph$. Then $(h,\ph)\m h*\ph$ is an $H$-module structure on $H^*$. 
Define an $\ca$-linear isomorphism $\L:H^*@>>>H$ by $\tc_w\m c_{ww_0}$. We show that second equality in (b) 
follows from the statement below:

(c) {\it $\L$ is $H$-linear.}
\nl
Let $h=T_{w_0}\i\in H$. Then ${}^th^\da=(-v^2)^{-l(w_0)}T_{w_0}$. By (c) we have 
$h*\tc_x=\sum_ya_{xw_0,yw_0}\tc_y$. Hence 
$$(-1)^{l(z)}\a_{xw_0,zw_0}=\sum_ya_{xw_0,yw_0}\tc_y(c_z)=\tc_x({}^th^\da c_z)=\tc_x(\sum_yb_{z,y}c_y)
=(-1)^{l(x)}b_{z,x}$$
and the second equality in (b) follows.

In the remainder of this subsection we prove (c). (This is a $\qq$-analogue of a result I proved in 1980 which is
reproduced in \cite{\BV, 2.25}.) It is enough to check that $\L(T_s*\tc_w)=T_s\L(\tc_w)$ for $w\in W,s\in\II$. 
Recall \cite{\KL} that there exists a symmetric function $\WW\T\WW@>>>\NN$, $y,w\m\mu(y,w)$ such that $\mu(y,w)=0$
unless $(-1)^{l(y)-l(w)}=-1$ and such that

(d) $T_sc_w=-c_w+\sum_{y;sy<y}\mu(y,w)vc_y$ if $sw>w$; $T_sc_w=v^2c_w$ if $sw<w$.
\nl
It follows that

$T_s\tc_w=-\tc_w$ if $sw>w$; 

$T_s\tc_w=v^2\tc_w+\sum_{y;sy>y}\mu(w,y)(-1)^{l(y)-l(w)}v\tc_y$ if $sw<w$.
\nl
The last equality can be written in the form 
$T_s\tc_w=v^2\tc_w-\sum_{y;sy>y}\mu(w,y)v\tc_y$. Since $T_s^\da=-v^2T_s\i$, we have

$T_s*\tc_w=-\tc_w+\sum_{y;sy>y}\mu(w,y)v\tc_y$ if $sw<w$; $T_s*\tc_w=v^2\tc_w$ if $sw>w$.
\nl
On the other hand, using (d) we have:

$T_sc_{ww_0}=-c_{ww_0}+\sum_{y;syw_0<yw_0}\mu(w_0y,w_0w)vc_{yw_0}$ if $sww_0>ww_0$, 

$T_sc_{ww_0}=v^2c_{ww_0}$ if $sww_0<ww_0$.
\nl
Now the condition that $sww_0>ww_0$ is equivalent to the condition that $sw<w$. Moreover, by \cite{\KL}, we have 
$\mu(y,w)=\mu(ww_0,yw_0)$ for any $y,w$. Hence $\L(T_s*\tc_w)=T_s\L(\tc_w)$. This proves (c) hence also (b).
 
\subhead 1.5\endsubhead
Let $\cd^{un}(Z)$ (resp. $\cd^{un}(Z')$) be the subcategory of $\cd(Z)$ (resp. $\cd(Z')$) whose objects are those
$L\in\cd(Z)$ (resp. $L\in\cd(Z')$) such that for any $i\in\ZZ$, any composition factor of $H^i(L)$ is isomorphic to
$L^w[[M_w]]$ (resp. $L'{}^w[[M_w]]$) for some $w\in\WW$. Let $\ck^{un}_\ca(Z)$ be the 
$\ca$-submodule of $\ck_\ca(Z)$ spanned by the basis elements $L^w[[M_w]]$. Let $\ck^{un}_\ca(Z')$ be the 
$\ca$-submodule of $\ck_\ca(Z')$ spanned by the basis elements $L'{}^w[[M_w]]$. We define an $\ca$-linear 
isomorphism $\Ps:H@>>>\ck^{un}_\ca(Z)$ by $(-1)^{-l(w)}c_w\m L^w[[M_w]]$ for all $w\in\WW$. We define an 
$\ca$-linear isomorphism $\Ps':H@>>>\ck^{un}_\ca(Z')$ by $(-1)^{-l(w)}c_w\m L'{}^w[[M_w]]$ for all $w\in\WW$. Let
$\t=\ft_!\fs^*:\cd(Z)@>>>\cd(Z')$, $\tit=\ft_*\fs^!:\cd(Z)@>>>\cd(Z')$. From the definitions we see that $\t$ 
restricts to a functor $\cd^{un}(Z)@>>>\cd^{un}(Z')$ denoted again by $\t$. Also $\cd_{un}(Z)$, $\cd^{un}(Z')$ are
stable under $\fD_Z,\fD_{Z'}$ hence $\tit$ restricts to a functor $\cd^{un}(Z)@>>>\cd^{un}(Z')$ denoted again by 
$\tit$.

Now if $L\in\cd^{un}(Z)\cap\cd_m(Z)$ then $\t(L),\tit(L)$ are naturally objects of $\cd^{un}(Z')\cap\cd_m(Z')$ and
$gr(L)\in\ck^{un}_\ca(Z)$, $gr(\t(L))\in\ck^{un}_\ca(Z')$, $gr(\tit(L))\in\ck^{un}_\ca(Z')$ are defined. Moreover,
there are well defined $\ca$-linear maps $gr(\t):\ck^{un}_\ca(Z)@>>>\ck^{un}_\ca(Z')$, 
$gr(\tit):\ck^{un}_\ca(Z)@>>>\ck^{un}_\ca(Z')$ such that $gr(\t)(gr(L))=gr(\t(L))$, $gr(\tit)(gr(L))=gr(\tit(L))$
for any $L$ as above. We show:

(a) {\it $gr(\t)(\Ps(h))=\Ps'(T_{w_0}h)$ for any $h\in H$.}
\nl
Let $\bfZ=\{(B_1,B,B')\in\cb^3;\po(B_1,B)=w_0\}$. We have a diagram with cartesian squares
$$\CD
Z@<\fs<<\fZ@>\ft>>Z'\\
@Vp_1VV     @Vp_2VV     @Vp_3VV\\
\cb^2@<\bfs<<\bfZ@>\bft>>\cb^2
\endCD$$
where $\bfs(B_1,B,B')=(B,B')$, $\bft(B_1,B,B')=(B_1,B')$ and $p_1,p_2,p_3$ are the obvious projections. Let 
$\cd^{un}(\cb^2)$ be the subcategory of $\cd(\cb^2)$ whose objects are those $L\in\cd(\cb^2)$ such that for any 
$i\in\ZZ$, any composition factor of $H^i(L)$ is equivariant for the diagonal $G^0$-action on $\cb^2$. Let 
$\ck^{un}_\ca(\cb^2)$ be the $\ca$-submodule of $\ck_\ca(\cb^2)$ spanned by the simple $G^0$-equivariant perverse
sheaves on $\cb^2$. Let $\bar\t=\bft_!\bfs^*:\cd^{un}(\cb^2)@>>>\cd^{un}(\cb^2)$. We define
$gr(\bar\t):\ck^{un}_\ca(\cb^2)@>>>\ck^{un}_\ca(\cb^2)$ in terms of $\bar\t$ in the same way as $gr(\t)$ was 
defined in terms of $\t$. Note that $p_1^*$, $p_3^*$ induce isomorphisms 
$\ck^{un}_\ca(\cb^2)@>\si>>\ck^{un}_\ca(Z)$, $\ck^{un}_\ca(\cb^2)@>\si>>\ck^{un}_\ca(Z')$. Moreover, from the 
cartesian diagram above we see that $p_3^*\bar\t=\t p_1^*$. We see that (a) is reduced to the well known 
description of $gr(\bar\t)$ in terms of left multiplication by $T_{w_0}$ in $H$.

We show:

(b) {\it $gr(\tit)(\Ps(h))=\Ps'(T_{w_0}\i h)$ for any $h\in H$.}
\nl
Note that $\fD_Z,\fD_{Z'}$ induce involutions of $\ck^{un}_\ca(Z)$, $\ck^{un}_\ca(Z')$ which are semilinear with 
respect to the ring involution $v^i\m v^{-i}$ of $\ca$ and are denoted again by $\fD_Z,\fD_{Z'}$. We have 
$\fD_Z\Ps(h)=\Ps(\bar h)$, $\fD_{Z'}\Ps'(h)=\Ps'(\bar h)$ for all $h\in H$. Moreover we have
$\tit=\fD_{Z'}\t\fD_Z$. Hence $gr(\tit)=\fD_{Z'}gr(\t)\fD_Z$. Using (a) we see that
$$\align&gr(\tit)(\Ps(h))=\fD_{Z'}gr(\t)\fD_Z(\Ps(h))=\fD_{Z'}gr(\t)(\Ps(\bar h))
\\&=\fD_{Z'}\Ps'(T_{w_0}\bar h)=\Ps'(\ov{T_{w_0}\bar h})=\Ps'(T_{w_0}\i h),\endalign$$
as required.

\subhead 1.6\endsubhead
Let $\l_1,\l_2,\do,\l_k$ be algebraic numbers in $\bbq$ such that for $i\in[1,k]$ any complex conjugate of 
$\l_i$ has absolute value $q^{t_i/2}$ where $t_i\in\NN$ and let $e_1,e_2,\do,e_k$ in $\{1,-1\}$ be such that 
$e_1\l_1^s+e_2\l_2^s+\do+e_k\l_k^s=0$ for all $s\in\ZZ_{>0}$. Then we have 
$e_1v^{t_1}+e_2v^{t_2}+\do+e_kv^{t_k}=0$ in $\ca$. The proof is left to the reader.

\subhead 1.7\endsubhead
Recall the assumption 1.2(a). Let $s\in\ZZ_{>0}$. By results in \cite{\CS} (when $D=G^0$) and \cite{\CDG, X} (in 
the general case) for any $A,A'\in\Xi$ we have the orthogonality relation:
$$\sum_{\x\in D_s}\c_{\uA,s}(\x)\c_{\uA',s}(\x)=\d|G^0_s|\o_A^s\tag a$$
where $\d=1$ if $A'=A^*$ and $\d=0$, otherwise; $\o_A$ is a root of $1$ in $\bbq$ not depending on $s$; moreover, 

(b) {\it if $A\in\Xi$, $w\in\WW$ and $j\in\ZZ$ are such that $(A:H^j(\bK^w_D))>0$ then $j=d_A\mod2$.}
\nl
By the relative hard Lefschetz theorem of Deligne \cite{\BBD, 5.4.10} applied to the projective morphism $\fd$ 
and to the simple perverse sheaf $\fc^*L^w[\D+l(w)]$ on $\tZ$, we see that for any $w\in\WW$, $j\in\ZZ$ we have

(c) $H^j(\bK^w_D)\cong H^{2\D+2\l(w)-j}(\bK^w_D)$.
\nl
From the fact that $\fd_!$ commutes with Verdier duality we see that $\fD(\bK^w_D)\cong\bK^w_D[2\D+2l(w)]$ hence
$\fD(H^j(\bK^w_D)\cong H^{2\D+2l(w)-j}(\bK^w_D)$. Combining this with (c) we obtain
$\fD(H^j(\bK^w_D)\cong H^j(\bK^w_D)$. Hence for any $A\in\Xi$ we have 

(d) $(\fD_D(A):H^j(\bK^w_D))=(A:H^j(\bK^w_D))$.

\subhead 1.8\endsubhead
Let $\cd^{un}(D)$ be the subcategory of $\cd(D)$ whose objects are those $K\in\cd(D)$ such that for any $i\in\ZZ$,
any composition factor of $H^i(K)$ is in $\hD^{un}$. For $w\in\WW$, $\bbq^w,L^w$ come naturally from objects 
$\ubbq^w,\uL^w$ of $\cd_m(Z)$ such that Frobenius acts trivially on $\ch^0(\ubbq^w)_x$, $\ch^0(\uL^w)_x$ for any 
$\FF_q$-rational point $x$ of $Z^w$. Then we have naturally $\uL^w[[M_w]]\in\cd_m(Z)$ (it is perverse, pure of 
weight $0$). We set $\ubK^w_D=\fd_!\fc^*(\uL^w)\in\cd_m(D)$ (it is pure of weight zero) and 
$\uK^w_D=\fd_!\fc^*(\ubbq^w)\in\cd_m(D)$, so that $gr(\ubK^w_D)\in\ck^{un}_\ca(D)$, 
$gr(\uK^w_D)\in\ck^{un}_\ca(D)$ are defined. 

Let $s\in\ZZ_{>0}$. From \cite{\KLL} we can deduce
$$\c_{\uL^w,s}=\sum_{y\in\WW}P_{y,w}(q^s)\c_{\ubbq^y,s}.\tag a$$
Applying the linear map $\fd_{!s}\fc^*_s$ to both sides we obtain
$$\fd_{!s}\fc^*_s(\c_{\uL^w,s})=\sum_{y\in\WW}P_{y,w}(q^s)\fd_{!s}\fc^*_s(\c_{\ubbq^y,s}).$$
We have 
$$\c_{\uK^w_D,s}=\fd_{!s}\fc^*_s(\c_{\ubbq^w,s}),\c_{\ubK^w_D,s}=\fd_{!s}\fc^*_s(\c_{\uL^w,s}).$$ 
It follows that
$$\c_{\ubK^w_D,s}=\sum_{y\in\WW}P_{y,w}(q^s)\c_{\uK^y_D,s}.\tag b$$
Let $A\in\Xi$. For any $j\in\ZZ$ the mixed perverse sheaf $H^j(\ubK^w_D)$ (pure of weight $j$) is canonically of 
the form $\op_{A\in\Xi}V'_{j,A}\ot\uA$ where $V'_{j,A}$ are finite dimensional vector spaces on which the 
Frobenius map acts naturally with a (multi)set of eigenvalues $E'_{j,A}$ in $\bbq$ such that each 
$\l\in E'_{j,A}$ is an algebraic number all of whose complex conjugates have absolute value $q^{(j-\De)/2}$. For 
any $s\in\ZZ_{>0}$ and any $j\in\ZZ$ we have
$$\c_{H^j(\ubK^w_D),s}=\sum_{A\in\Xi}\sum_{\l\in E'_{j,A}}\l^s\c_{\uA,s}$$
as functions on $D_s$. It follows that
$$\c_{\ubK^w_D,s}=\sum_{A\in\Xi}\ct'_{w,A;s}\c_{\uA,s}$$
where 
$$\ct'_{w,A;s}=\sum_{j\in\ZZ}(-1)^j\sum_{\l\in E'_{j,A}}\l^s.$$
We set
$$\ct'_{w,A}=\sum_{j\in\ZZ}(-1)^j\dim V'_{j,A}v^j=\sum_{j\in\ZZ}(-1)^j(A:H^j(\bK^w_D))v^{j-\De}.
$$
Note that
$$\ct'_{w,A}=v^{-\De}(A:gr(\ubK^w_D)).$$
For any $i,j\in\ZZ$ the mixed perverse sheaf $H^i(\uK^w_D)_j$ (pure of weight $j$) is canonically of the form
$\op_{A\in\Xi}V''_{i,j,A}\ot\uA$ where $V''_{i,j,A}$ are finite dimensional vector spaces on which the Frobenius 
map acts naturally with a (multi)set of eigenvalues $E''_{i,j,A}$ in $\bbq$ such that each $\l\in E''_{i,j,A}$ is
an algebraic number all of whose complex conjugates have absolute value $q^{(j-\De)/2}$. For any $s\in\ZZ_{>0}$ 
and any $i,j\in\ZZ$ we have
$$\c_{H^i(\uK^w_D)_j,s}=\sum_{A\in\Xi}\sum_{\l\in E''_{i,j,A}}\l^s\c_{\uA,s}$$
as functions on $D_s$. It follows that
$$\c_{\uK^w_D,s}=\sum_{A\in\Xi}\ct''_{w,A;s}\c_{\uA,s}$$
where 
$$\ct''_{w,A;s}=\sum_{i,j\in\ZZ}(-1)^i\sum_{\l\in E''_{j,A}}\l^s.$$
We set
$$\ct''_{w,A}=\sum_{i,j\in\ZZ}(-1)^i\dim V''_{j,A}v^j
\sum_{i,j\in\ZZ}(-1)^i(A:H^i(\uK^w_D)_j)v^{j-\De}.$$
Note that
$$\ct''_{w,A}=v^{-\De}(A:gr(\uK^w_D)).$$
Using now (b) we see that
$$\sum_{A\in\Xi}\ct'_{w,A;s}\c_{\uA,s}=
\sum_{y\in\WW}P_{y,w}(q^s)\sum_{A\in\Xi}\ct''_{w,A;s}\c_{\uA,s}.$$
Since the functions $\c_{\uA,s}$ (with $A\in\Xi$) are linearly independent (see 1.7(a)), it follows that
$$\ct'_{w,A;s}=\sum_{y\in\WW}P_{y,w}(q^s)\ct''_{y,A;s}\tag c$$
for any $w\in\WW,A\in\Xi$.
Applying 1.6 to (c) we obtain the equality
$$\ct'_{w,A}=\sum_{y\in\WW}P_{y,w}(v^2)\ct''_{y,A}\tag d$$
in $\ca$. Equivalently,
$$(A:gr(\ubK^w_D))=\sum_{y\in\WW}P_{y,w}(v^2)(A:gr(\uK^w_D)).\tag e$$
We define an $\ca$-linear map $\Ph:H@>>>\ck^{un}_\ca(D)$ by $\Ph(T_w)=gr(\uK^w_D)$ for all
$w\in\WW$. Now (e) shows that 
$$\Ph(v^{l(w)}c_w)=gr(\ubK^w_D)\tag f$$ 
for all $w\in\WW$.

\subhead 1.9 \endsubhead
We show that for $x,y$ in $W$ we have
$$\sum_{\z\in Z_s}\c_{\uL^x,s}(\z)\c_{\uL^y,s}(\z)
=|G^0_s|\sum_{z\in\WW}P_{z,x}(q^s)P_{z,y}(q^s)q^{-l(w_0z)s}.\tag a$$
Using 1.8(a) we see that the left hand side equals
$$\align&\sum_{x',y'\in\WW}P_{x',x}(q^s)P_{y',y}(q^s)\sum_{\z\in Z_s}
\c_{\ubbq^{x'},s}(\z)\c_{\ubbq^{y'},s}(\z)\\&=
\sum_{x',y'\in\WW}P_{x',x}(q^s)P_{y',y}(q^s)|\{\z\in Z_s\cap Z^{x'}\cap Z^{y'}\}|
=\sum_{z\in\WW}P_{z,x}(q^s)P_{z,y}(q^s)|Z_s^z|\endalign$$
and it remains to use the equality $|Z_s^z|=|G^0_s|q^{-l(w_0z)s}$ for any $z\in\WW$.

\subhead 1.10\endsubhead
From \cite{\PAR, 6.5} we see that $\fe:\fc_!\fd_*:\cd(D)@>>>\cd(Z)$ restricts to a functor
$\cd^{un}(D)@>>>\cd^{un}(Z)$ denoted again by $\fe$. Let $A\in\Xi$. We have $\fe(A)\in\cd^{un}(Z)$ and 
$\fe(\uA)\in\cd_m(Z)$ hence $gr(\fe(\uA))\in\ck^{un}_\ca(Z)$ is defined. We have the following result:
$$gr(\fe(\uA))=\Ps((-1)^{d_A}T_{w_0}\sum_{x\in\WW}
v^{-l(w_0x)}(A^\bul:gr(\ubK^{w_0x}_D)))(-1)^{l(x)}c_x).\tag a$$
For any $i,j\in\ZZ$ the mixed perverse sheaf $H^i(\fe(\uA))_j$ (pure of weight $j$) is canonically of the form
$\op_{x\in\WW}V_{i,j,x}\ot\uL^x[[M_x]]$ where $V_{i,j,x}$ are finite dimensional vector spaces on which the 
Frobenius map acts naturally with a (multi)set of eigenvalues $E_{i,j,x}$ in $\bbq$ such that each 
$\l\in E_{i,j,x}$ is an algebraic number all of whose complex conjugates have absolute value $q^{j/2}$. For any 
$s\in\ZZ_{>0}$ and any $i,j\in\ZZ$ we have
$$\c_{H^i(\fe(\uA))_j,s}=\sum_{x\in\WW}\sum_{\l\in E_{i,j,x}}\l^s\c_{\uL^x[[M_x]],s}$$
as functions on $Z_s$. It follows that
$$\c_{\fe(\uA),s}=\sum_{x\in\WW}S_{x,A;s}\c_{\uL^x[[M_x]],s}$$
where
$$S_{x,A;s}=\sum_{i,j\in\ZZ}(-1)^i\sum_{\l\in E_{i,j,x}}\l^s.$$
We set
$$S_{x,A}=\sum_{i,j\in\ZZ}(-1)^i\dim V_{i,j,x}v^j
=\sum_{i,j\in\ZZ}(-1)^i(\bbq^{x\sh}[[M_x]]:H^i(\fe(\uA))_j)v^j.$$
Note that
$$S_{x,A}=(L^x[[M_x]]:gr(\fe(\uA))).$$
For $w\in\WW$ we have
$$\align&\sum_{\x\in D_s}\c_{\uA,s}(\x)\c_{\ubK^w,s}(\x)=
\sum_{\x\in D_s}\c_{\uA,s}(\x)\sum_{\et\in\tZ_s;\fd(\et)=\x}\c_{\fc^*\uL^w,s}(\et)\\&
=\sum_{\et\in\tZ_s}\c_{\uA,s}(\fd(\et))\c_{\uL^w,s}(\fc(\et)),\endalign$$
$$\align&\sum_{\z\in Z_s}\c_{\fc_!\fd^*\uA,s}(\z)\c_{\uL^w,s}(\z)\\&=
\sum_{\z\in Z_s}\sum_{\et\in\tZ_s;\fc(\et)=\z}\c_{\fd^*\uA,s}(\et)\c_{\uL^w,s}(\z)=
\sum_{\et\in\tZ_s}\c_{\fd^*\uA,s}(\et)\c_{\uL^w,s}(\fc(\et)).\endalign$$
It follows that
$$\sum_{\x\in D_s}\c_{\uA,s}(\x)\c_{\ubK^w,s}(\x)=\sum_{\z\in Z_s}\c_{\fc_!\fd^*\uA,s}(\z)\c_{\uL^w,s}(\z).$$
We rewrite this as follows
$$\align&\sum_{\x\in D_s}\c_{\uA,s}(\x)\sum_{A'\in\Xi}\ct'_{w,A';s}\c_{\uA',s}(\x)=
\sum_{\z\in Z_s}\sum_{x\in\WW}S_{x,A;s}\c_{\uL^x[[M_x]],s}(\z)\c_{\uL^w,s}(\z)\\&
=\sum_{x\in\WW}S_{x,A;s}(-q^{s/2})^{-M_x}\sum_{z\in\WW}P_{z,x}(q^s)P_{z,w}(q^s)|G^0_s|q^{-l(w_0z)s}\endalign$$
where the last equality follows from 1.9(a) and the equality
$$\c_{\uL^x[[M_x]],s}=(-q^{s/2})^{-M_x}\c_{\uL^x,s}.$$
Using now 1.7(a) we deduce
$$|G^0_s|\ct'_{w,A^*;s}\o_A^s=
|G^0_s|\sum_{x,z\in\WW}S_{x,A;s}(-q^{s/2})^{-M_x}P_{z,x}(q^s)P_{z,w}(q^s)q^{-l(w_0z)s}$$
or equivalently (see 1.8(c)):
$$\sum_{z\in\WW}P_{z,w}(q^s)\ct''_{z,A^*;s}\o_A^s
=\sum_{x,z\in\WW}S_{x,A;s}(-q^{s/2})^{-M_x}P_{z,x}(q^s)P_{z,w}(q^s)q^{-l(w_0z)s}.$$
We multiply both sides by $Q_{w,u}(q^s)$ (entries of the inverse matrix of $(P_{y,w}(q^s))$
and sum over $w\in\WW$. We obtain
$$\ct''_{u,A^*;s}\o_A^s=\sum_{x\in\WW}S_{x,A;s}(-q^{s/2})^{-M_x}P_{u,x}(q^s)q^{-l(w_0u)s}.$$
We multiply both sides by $Q_{y,u}(q^s)q^{l(w_0u)s}$ and sum over $u\in\WW$. We obtain
$$S_{y,A;s}(-q^{s/2})^{-M_y}=\sum_{u\in\WW}Q_{y,u}(q^s)q^{l(w_0u)s}\ct''_{u,A^*;s}\o_A^s.$$
Applying 1.6 to the previous equality we obtain
$$S_{y,A}(-v)^{-M_y}=\sum_{u\in\WW}Q_{y,u}(v^2)v^{2l(w_0u)}\ct''_{u,A^*}.$$
Here we substitute $\ct''_{u,A^*}=v^{-\De}(A^*:gr(\ubK^u_D)))=v^{-\De}(A^*:\Ph(T_u))$ and 
$Q_{y,u}=(-1)^{l(y)-l(u)}P_{w_0u,w_0y}$. Note also that 
$(A^*:gr(\ubK^u_D)))=(A:gr(\ubK^u_D)))$, by 1.7(d). We obtain
$$\align&S_{y,A}=v^{-\De}(-v)^{M_y}(A:\Ph(\sum_{u\in\WW}(-1)^{l(y)-l(u)}P_{w_0u,w_0y}(v^2)v^{2l(w_0u)}T_u))\\&
=v^{-\De}(-v)^{M_y}(\dd(A):\dd(\Ph(\sum_{u\in\WW}(-1)^{l(y)-l(u)}P_{w_0u,w_0y}(v^2)v^{2l(w_0u)}T_u))).\endalign$$
(We use that $(\r:\r')=(\dd(\r):\dd(\r'))$ for any $\r,\r'$ in $\ck^{un}_\ca(D)$.) From \cite{\CDG, IX, 42.9} we 
have 
$\dd(\Ph(h))=\Ph(h^\da)$ for any $h\in H$. Hence the previous formula for $S_{y,A}$ becomes 
$$\align&S_{y,A}=v^{-\De}(-v)^{M_y}(-1)^{d'_A}(-1)^{l(y)}\\&\T
(A^\bul:\Ph(\sum_{u\in\WW}(-1)^{l(u)}P_{w_0u,w_0y}(v^2)v^{2l(w_0u)}(-v^2)^{l(u)}T_{u\i}\i)).\endalign$$
We now replace 

$T_{u\i}\i$ by $T_{w_0}\i T_{w_0u}$, 

$(-1)^{l(u)}v^{2l(w_0u)}(-v^2)^{l(u)}$ by $v^{2l(w_0)}$,

$\sum_{u\in\WW}P_{w_0u,w_0y}(v^2)T_{w_0u}$ by $v^{l(w_0y)}c_{w_0y}$,

$v^{-\De}(-v)^{M_y}(-1)^{d'_A}(-1)^{l(y)}$ by $(-1)^{l(w_0)}v^{-\l(w_0y)}(-1)^{d_A}$; 
\nl
we obtain
$$S_{y,A}=(-v^2)^{l(w_0)}(-1)^{d_A}(A^\bul:\Ph(T_{w_0}\i c_{w_0y})).$$
We have
$$gr(\fe(\uA))=\sum_{y\in\WW}S_{y,A}L^y[[M_y]]$$
Hence
$$\align&\Ps\i(gr(\fe(\uA)))=\sum_{y\in\WW}S_{y,A}(-1)^{l(y)}c_y\\&
=(-v^2)^{l(w_0)}(-1)^{d_A}\sum_{y\in\WW}(A^\bul:\Ph(T_{w_0}\i c_{w_0y}))(-1)^{l(y)}c_y.\endalign$$
Using now 1.4(a), 1.4(b) we obtain
$$\align&\Ps\i(gr(\fe(\uA)))
=(-v^2)^{l(w_0)}(-1)^{d_A}\sum_{x,y\in\WW}(A^\bul:\Ph(a_{w_0y,w_0x}c_{w_0x}))(-1)^{l(y)}c_y\\&
=(-v^2)^{l(w_0)}(-1)^{d_A}\sum_{x,y\in\WW}(A^\bul:\Ph((-1)^{l(y)-l(x)}b_{x,y}c_{w_0x}))(-1)^{l(y)}c_y.\endalign$$
Here we substitute $(-v^2)^{l(w_0)}\sum_{y\in\WW}b_{x,y}c_y=T_{w_0}c_x$ (see 1.4(a)); we obtain
$$\Ps\i(gr(\fe(\uA)))=(-1)^{d_A}\sum_{x\in\WW}(A^\bul:\Ph((-1)^{l(x)}c_{w_0x}))T_{w_0}c_x.$$
We now use 1.8(f) and apply $\Ps$ to both sides; we obtain (a).

\subhead 11. Proof of 1.2(b)\endsubhead
Let $A\in\Xi$. From the definitions we have 
$$gr(\b(\uA))=gr(\tit\fe(\uA))=gr(\tit)(gr(\fe(\uA))).$$
We write the equation 1.5(b) for $h=\Ps\i(gr(\fe(\uA)))\in H$. We obtain 
$$gr(\tit)(gr(\fe(\uA)))=\Ps'(T_{w_0}\i \Ps\i(gr(\fe(\uA)))).$$ Using now 1.10(a) we obtain
$$\align&gr(\b(uA))=\Ps'((-1)^{d_A}\sum_{x\in\WW}v^{-l(w_0x)}(A^\bul:gr(\ubK^{w_0x}_D)))(-1)^{l(x)}c_x)\\&
=\sum_{x\in\WW}\sum_{j\in\ZZ}(-1)^{j+d_A}(A^\bul:H^j(\bK^{w_0x}_D))v^{j-l(w_0x)}L'{}^x[[M_x]].\endalign$$
By 1.7(b) we have $(A^\bul:H^j(\bK^{w_0x}_D))=0$ unless $j+d_A=0\mod2$ (note that $d_A=d_{A^\bul}$). Hence we have
$$gr(\b(uA))=\sum_{x\in\WW}\sum_{j\in\ZZ}(A^\bul:H^j(\bK^{w_0x}_D))v^{j-l(w_0x)}L'{}^x[[M_x]].\tag a$$
Now 1.2(b) follows.

\subhead 1.12. The functor $\ti\b$\endsubhead
Let $\ti\b=\ft_!\fs^*\fc_*\fd^!:\cd(D)@>>>\cd(Z')$. Let $A\in\Xi$. The following equality suggests that $\ti\b$ 
might be equal to $\b$ up to a twist:
$$gr(\ti\b(\uA))=v^{-2\De}gr(\b(\uA)).\tag a$$
Note that $\ti\b(A)=\fD_{Z'}\b(\fD_DA)$. It is enough to show that
$$gr(\fD_{Z'}\b(\uA^*))=v^{-2\De}gr(\b(\uA))$$
or, by 1.11(a), that
$$\align&\sum_{x\in\WW}\sum_{j\in\ZZ}((\fD_DA)^\bul:H^j(\bK^{w_0x}_D))v^{-j+l(w_0x)}L'{}^x[[M_x]]\\&
=\sum_{x\in\WW}\sum_{j\in\ZZ}(A^\bul:H^j(\bK^{w_0x}_D))v^{-2\De+j-l(w_0x)}L'{}^x[[M_x]].\endalign$$
We have $(\fD_DA)^\bul=\fD_D(A^\bul)$. Hence by 1.7(c), 1.7(d) we have for any $x\in\WW$:
$$((\fD_DA)^\bul:H^j(\bK^{w_0x}_D))=(A^\bul:H^{2\De+2l(w_0x)-j}(\bK^{w_0x}_D)).$$
Hence it is enough to observe that for any $x\in\WW$ we have
$$\sum_{j\in\ZZ}(A^\bul:H^{2\De+2l(w_0x)-j}(\bK^{w_0x}_D))v^{-j+l(w_0x)}
=\sum_{j'\in\ZZ}(A^\bul:H^{j'}(\bK^{w_0x}_D))v^{-2\De+j'-l(w_0x)}.$$
(We use the substitution $j'=2\De+2l(w_0x)-j$.)

\subhead 1.13\endsubhead
To simplify the notation, in the remainder of \S1 we assume that $D=G^0$. (Similar results hold without this 
assumption.) Let $\ce$ be a set of representatives for the isomorphism classes of simple $\QQ[\WW]$-modules; for 
each $E\in\ce$ let $E^v$ be the corresponding simple $H^v$-module where $H^v=\QQ(v)\ot_\ca H$. Let $A\in\Xi$.
By \cite{\CS, III, 14.11}, for any $w\in\WW$ we have
$$\sum_{j\in\ZZ}(A:H^j(\bK^w_D))v^j=(-1)^{d'_A}\sum_{E\in\ce}v^{\De+l(w)}\g_{A,E}\tr(c_w,E^v)$$
where $\g_{A,E}$ are certain rational numbers. Taking this into account we can rewrite 1.11(a) as follows:
$$\Ps'{}\i(gr(\b(uA)))=(-1)^{d'_A}v^\De\sum_{E\in\ce}^r\g_{A^\bul,E}\fC_E$$
(equality in $H^v$) where
$$\fC_E=\sum_{x\in\WW}(-1)^{l(x)}\tr(c_{w_0x},E^v)c_x\in H.\tag a$$
It is known that $(-1)^{d'_A}\g_{A^\bul,E}=\g_{A,E^\da}$ where $E^\da\in\ce$ is isomorphic to $E\ot\sgn$ and $\sgn$
is the sign representation of $\WW$. Hence we have
$$\Ps'{}\i(gr(\b(uA)))=v^\De\sum_{E\in\ce}^r\g_{A,E^\da}\fC_E.\tag b$$
Note that in the sum over $E$ in (b) can be restricted to the $E$ which belong to a fixed two-sided cell 
(depending on $A$); this is a known property of the coefficients $\g_{A,E^\da}$. We show that for any $E\in\ce$ we
have
$$T_{w_0}\i\fC_E\in\text{centre}(H), T_{w_0}\fC_E\in\text{centre}(H).\tag c$$
(These two statements are equivalent since $T_{w_0}^2$ is in the centre of $H$.) It is enough to show that the 
image of $(T_{w_0}\i\fC_E)^\da$ is in the centre of $H$. We have
$$T_{w_0}\i\fC_E
=\sum_{x,y,z\in\WW}(-1)^{l(x)}v^{-l(w_0x)}v^{-l(x)}T_{w_0}\i P_{w_0y,w_0x}P_{y,x}\tr(T_{w_0y},E^v)T_z.$$
Using 1.3(a) we obtain
$$\align&T_{w_0}\i\fC_E=\sum_{y,z\in\WW}(-1)^{l(y)}v^{-l(w_0)}T_{w_0}\i\d_{y,z}\tr(T_{w_0y},E^v)T_z\\&
=\sum_{y\in\WW}(-1)^{l(y)}v^{-l(w_0)}\tr(T_{w_0y},E^v)T_{w_0}\i T_y
=\sum_{y\in\WW}(-1)^{l(y)}v^{-l(w_0)}\tr(T_{w_0y},E^v)T_{y\i w_0}\i\\&
=\sum_{u\in\WW}(-1)^{l(w_0u)}v^{-l(w_0)}\tr(T_u,E^v)T_{u\i}\i.\endalign$$
Hence
$$(T_{w_0}\i\fC_E)^\da=(-v)^{-l(w_0)}\fC'_E$$
where
$$\fC'_E=\sum_{u\in\WW}v^{-2l(u)}\tr(T_{u\i},E^v)T_u.$$
It is well known that the elements $\fC'_E (E\in\ce)$ for a basis of the centre of $H^v$. This proves (c).

\subhead 1.14\endsubhead
If $X$ is an algebraic variety and $K\in\cd(X)$ we set 
$$gr_1(K)=\sum_{i\in\ZZ}(-1)^iH^i(K)\in\ck(X).$$
 Assume that $A\in\Xi$ is cuspidal. From 1.11(a) we have
$$gr_1(\b(A))=(-1)^{d_A}\sum_{x\in\WW}(A:gr_1(\bK^{w_0x}_D))L'{}^x[[M_x]].\tag a$$
(We have $A^\bul=A$ since $A$ is cuspidal.) Let $\G_A$ be the set of all $x\in\WW$ such that $L'{}^x[[M_x]]$ 
appears with $\ne0$ coefficient in the sum (a) and $x$ has maximum possible length with this property. Note that 
$\G_A\ne\em$. We show:

(b) {\it $w_0\G_A$ is contained in a single conjugacy class in $\WW$.}
\nl
Let $\G'_A$ be the set of all $w\in\WW$ such that $(A:gr_1(\bK^w_D)\ne0\}$ and $l(w)$ is minimum possible with 
this property. We have $w_0\G_A=\G'_A$. Let $\G''_A$ be the set of all $w\in\WW$ such that $(A:gr_1(K^w_D)\ne0\}$
and $l(w)$ is minimum possible with this property. Since $gr_1(\bK^w_D)=\sum_{y\in\WW;y\le w}P_{y,w}(1)gr(K^y_D)$,
we see that $\G'_A=\G''_A$. It is enough to show that $\G''_A$ is contained in a single conjugacy class in $\WW$.
This follows from the Corollary to Theorem 2.18 in \cite{\RA}. For a description of the conjugacy classes in $\WW$
that arise in this manner, see \cite{\RA}.

\head 2. Tensor products of character sheaves\endhead
\subhead 2.1\endsubhead
An element $g\in G$ is said to be {\it isolated} in $G$ if there is no proper parabolic subgroup $P$ of $G^0$ with
Levi $L$ such that $gPg\i=P,gLg\i=L$ and $L$ and $Z_G(g_s)^0\sub L$, see \cite{\CDG, I, 2.2}. A subset $C$ of $G$
is said to be an {\it isolated stratum} of $G$ if $C$ is contained in a connected component $D'$ of $G$, $C$ is a
single orbit of the action $(z,x):y\m xzyx\i$ of ${}^{D'}\cz^0_{G^0}\T G^0$ on $D'$ and if some/any element of $C$
is isolated in $G$, see \cite{\CDG, I, 3.3}. If $C$ is an isolated stratum of $G$ and $n\in\NN_\kk^*$, let 
$\cs_n(C)$ be the category whose objects are the local systems on $C$ that are equivariant for the transitive 
${}^{D'}\cz^0_{G^0}\T G^0$-action $(z,x):y\m xz^nyx\i$ on $D'$ ($D'$ is associated to $C$ as above). Let $\cs(C)$
be the category whose objects are the local systems on $C$ that are in $\cs_n(C)$ for some $n$ as above.

Following \cite{\CDG, I, 3.5}, let $\AA$ be the set of all pairs $(L,S)$ where $L$ is a Levi subgroup of some 
parabolic of $G^0$ and $S$ is an isolated stratum of $N_G(L)$ with the following property: there exists a 
parabolic subgroup $P$ of $G^0$ with Levi $L$ such that $S\sub N_GP$. For $(L,S)\in\AA$ let 
$\tcw_S=\{n\in N_{G^0}L;nSn\i=S\}$ (a subgroup of $N_{G^0}L$) and $\cw_S=\tcw_S/L$ (a subgroup of the finite group
$N_{G^0}L/L$). Now $\tcw_S$ acts on $S$ by conjugation. Hence if $\nu\in\tcw_S$ and $\ce\in\cs(S)$ then $\nu^*\ce$
is a well defined local system (necessarily in $\cs(S)$).

For $(L,S)\in\AA$ let $S^*=\{g\in S;Z_G(g_s)^0\sub L\}$ (an open dense subset of $S$, see \cite{\CDG, I, 3.11}) 
and let $Y_{L,S}=\cup_{x\in G^0}xS^*x\i$. By \cite{\CDG, I, 3.16}, $Y_{L,S}$ is a locally closed irreducible 
subvariety of $G$. Now $G^0$ acts on $\AA$ by conjugation; moreover $(L,S)\in\AA$ and $(L',S')\in\AA$ are in the 
same $G^0$-orbit if and only if $Y_{L,S}=Y_{L',S'}$, see \cite{\CDG, I, 3.12}. The subsets $Y_{L,S}$ with $(L,S)$
running through a set of representatives for the $G^0$-orbits in $\AA$ form a partition of $G$ into finitely many
subsets called the {\it strata} of $G$, see \cite{\CDG, I, 3.12}. By \cite{\CDG, I, 3.15}, the closure of any 
stratum of $G$ is a union of strata of $G$. For $(L,S)\in\AA$ let 
$\tY_{L,S}=\{(g,xL)\in G^0\T G/L;x\i gx\in S^*\}$. Define $\p:\tY_{L,S}@>>>Y_{L,S}$ by $(g,xL)\m g$. Now $\cw_S$ 
acts freely on $\tY_{L,S}$ by $\nu:(g,xL)\m(g,x\nu\i L)$. This makes $\p:\tY_{L,S}@>>>Y_{L,S}$ into a principal 
$\cw_S$-bundle, see \cite{\CDG, I, 3.13}. By \cite{\CDG, I, 3.17}, $Y_{L,S}$ and $\tY_{L,S}$ are smooth. Let 
$(L,S)\in\AA$. We have a diagram $\tY_{L,S}@<a<<R@>b>>S$ where $R=\{(g,x)\in G\T G^0;x\i gx\in S^*\}$ and 
$a(g,x)=(g,xL)$, $b(g,x)=x\i gx$. Let $\ce\in\cs(S)$. There is a well defined local system $\tce$ on $\tY_{L,S}$ 
such that $b^*\ce=a^*\tce$. Now $\tcw_S$ acts on $\tY_{L,s}$ through its quotient $\cw_S$, on $Z$ by
$\nu:(g,x)\m(g,x\nu\i)$ and on $S$ as above. These actions are compatible with $a,b$. Hence if $\nu\in\tcw_S$ 
represents $w\in\cw_S$ we have $w^*\tce\cong\wt{\nu^*\ce}$ where $w^*\tce$ is defined using the $\cw_S$-action on
$\tY_{L,S}$. Note that $\p_!\tce$ is a local system on $Y_{L,S}$. 

Let $Y$ be a stratum of $G$. Let $\Lo(Y)$ be the category whose objects are the local systems on $Y$ that are 
isomorphic to a direct summand of the local system $\p_!\tce$ for some $\ce\in\cs(S)$ (where 
$Y=Y_{L,S},(L,S)\in\AA$); this is independent of the choice of $(L,S)$ such that $Y=Y_{L,S}$. Note that any object
of $\Lo(Y)$ is semisimple. The local system $\bbq$ on $Y$ belongs to $\Lo(Y)$.
Clearly $\Lo(Y)$ is closed under direct sum. We show:

(a) {\it $\Lo(Y)$ is closed under $\ot$.}
\nl
Let $(L,S)$ be such that $Y=Y_{L,S}$. Let $\cl,\cl'$ be objects of $\Lo(Y)$. By assumption we can find 
$\ce\in\cs(S)$, $\ce'\in\cs(S)$ such that $\cl$ is a direct summand of $\p_!\tce$ and $\cl'$ is a direct summand 
of $\p_!\tce'$. Then $\cl\ot\cl'$ is a direct summand of the local system $\p_!\tce\ot\p_!\tce'$. Hence it is 
enough to show that $\p_!\tce\ot\p_!\tce'\in\Lo(Y)$. We have
$\p_!\tce\ot\p_!\tce'=\op_{w\in\cw_S}\p_!(\tce\ot w^*\tce')$. Hence it is enough to show that for any $w\in\cw_S$,
$\p_!(\tce\ot w^*\tce')\in\Lo(Y)$. Thus we must show that 
$\p_!(\tce\ot\wt{\nu^*\ce'})=\p_!(\wt{\ce\ot\nu^*\ce'})\in\Lo(Y)$ (where $\nu\in\tcw_S$ represents $w$); this 
holds since $\ce\ot\nu^*\ce'\in\cs(S)$. This proves (a).

\subhead 2.2\endsubhead
Let $Y$ be a stratum of $G$. Let $\fK(Y)$ be the Grothendieck group of the category $\Lo(Y)$. Note that any object
$\cl$ of $\Lo(Y)$ can be viewed as an element of $\fK(Y)$ denoted by $\ucl$. Let $\uLo(Y)$ be a set of 
representatives for the isomorphism classes of irreducible local systems in $\Lo(Y)$.

From 2.1(a) we see that $\fK(Y)$ is naturally a commutative ring in which for any $\cl,\cl'$ in $\Lo(Y)$ the 
product of $\ucl$, $\ucl'$ in $\fK(Y)$ is $\un{\cl\ot\cl'}$. This ring has a unit element: the class of the local system
$\bbq$ on $Y$.

Let $\fK_\ca(Y)=\ca\ot\fK(Y)$; this is naturally a commutative $\ca$-algebra with $1$. For any $\cl\in\Lo(Y)$ we
set $[\cl]=v^{\codim Y}\ucl\in\ck_\ca(Y)$. The elements $[\cl])$ where $\cl$ runs over $\uLo(Y)$ form an 
$\ca$-basis $B_Y$ of $\fK_\ca(Y)$.

\subhead 2.3\endsubhead
Let $\cd^*(G)$ be the subcategory of $\cd(G)$ whose objects are those $K\in\cd(G)$ such that for any stratum $Y$ 
of $G$ and any $i\in\ZZ$ we have $\ch^i(K)|_Y\in\Lo(Y)$. 

For any stratum $Y$ of $G$ and any irreducible local system $\cl$ in $\Lo(Y)$ the simple perverse sheaf 
$A_\cl:=IC(\bY,\cl)[\dim Y]^G$ is in $\cd^*(G)$, see \cite{\CDG, V, 25.2}. Conversely, let $A$ be a simple 
perverse sheaf in $\cd^*(G)$. We can find an irreducible local system $\cl'$ on an open dense smooth subvariety 
$V$ of $\supp(A)$ such that $A=IC(\supp(A),\cl')[\dim\supp(A)]^G$. The intersections of $\supp(A)$ with the 
various strata of $G$ form a partition of $\supp(A)$ into finitely many locally closed subvarieties; hence we can
find a stratum $Y$ of $G$ such that $Y\cap\supp(A)$ is open dense in $\supp(A)$. By assumption, 
$\ch^{-\dim\supp(A)}A|_Y\in\Lo(Y)$. Hence $\cl'|_{Y\cap V}=\cl_{Y\cap V}$. Since $Y\cap V$ is open dense in 
$\supp(A)$ we have $A=IC(\bY,\cl)[\dim Y]^G$. Moreover $\cl$ is automatically irreducible. We see that the simple
perverse sheaves in $\cd^*(G)$ are precisely the complexes of the form $A_\cl$ where $\cl$ is an irreducible local
system in $\Lo(Y)$ for some stratum $Y$ of $G$. 

\subhead 2.4\endsubhead
Let $\fK_\ca(G)=\op_Y\fK_\ca(Y)$; here $Y$ runs over the (finite) set of strata of $G$. We view $\fK_\ca(G)$ with
a commutative $\ca$-algebra structure which is the direct sum of the algebras $\fK_\ca(Y)$. Note that 
$B_G:=\sqc_YB_Y$ is an $\ca$-basis of $\fK_\ca(G)$.

Let $\cl\in\Lo(Y_0)$ where $Y_0$ is a stratum of $G$. Let $A=A_\cl$. If $s\in\ZZ_{>0}$ is sufficiently divisible 
then $F^{s*}A\cong A$ and we can choose an isomorphism $F^{s*}A@>\si>>A$ such that for any $g$ in an open dense 
subset of $\bY_0$ and any $s'\in\ZZ_{>0}s$ such that $F^{s'}(g)=g$, the eigenvalues of $F^{s'}$ on the stalk 
$\ch^{-dim\bY_0}(A)_g$ are roots of $1$ times $q^{s'\codim Y_0/2}$. We can also assume that each stratum of $G$ is
$F^s$-stable. For any stratum $Y$ of $G$ and any $i\in\ZZ$, $F^s$ induces an isomorphism 
$\ph_Y:F^{s*}\ch^i(A)|_Y@>\si>>\ch^i(A)|_Y$ and the local system $\ch^i(A)|_Y$ has a canonical filtration 
compatible with $\ph_Y$ whose subquotients $(\ch^i(A)|_Y)_j$ ($j\in\ZZ$) have the following property: for any 
$s'\in\ZZ_{>0}$ and any $g\in Y(F_{q^{s'}})$, any eigenvalue of $\ph_Y$ on the stalk $((\ch^i(A)|_Y)_j)_g$ is an
algebraic number all of whose complex conjugates have absolute value $q^{s'j/2}$. We set
$$\align&Gr(A)=(-1)^{\dim Y_0}\sum_Y\sum_{i,j\in\ZZ}(-1)^i\un{(\ch^i(A)|_Y)_j}v^j\\&
=(-1)^{\dim Y_0}\sum_Y\sum_{i,j\in\ZZ}(-1)^i[(\ch^i(A)|_Y)_j]v^{j-\codim Y}\in\fK_\ca(G)\endalign$$
where $Y$ runs over the strata of $G$. Note that $Gr(A)$ is independent of the choice of $s,\ph_Y$. Using Gabber's purity
theorem \cite{\BBD} and the fact that $\bY_0$ is a union of strata of $G$ (see \cite{\CDG, I, 3.15}) we see that 
in the second sum defininig $Gr(A)$, $j$ can be assumed to satisfy $j-\codim Y<0$ if $Y\sub\bY_0,Y\ne Y_0$ and 
$j-\codim Y=0$ if $Y=Y_0$. Thus $Gr(A_\cl)$ is equal to $[\cl]$ plus a $v\i\ZZ[v\i]$-linear combination of 
elements in $\sqc_{Y_1}B_{Y_1}$ where $Y_1$ runs over the strata of $G$ such that $Y_1\sub\bY_0$, $Y_1\ne Y_0$. We
see that the elements $Gr(A_\cl)$ with $[\cl]$ running through $B_G$ form an $\ca$-basis $\tB_G$ of $\fK_\ca(G)$.

Hence if $[\cl],[\cl']\in B_G$ we can write the product $Gr(A_\cl)Gr(A_{\cl'})$ in $\fK_\ca(G)$ uniquely in the 
form 
$$Gr(A_\cl)Gr(A_{\cl'})=\sum_{[\cl'']\in B_G}f_{\cl,\cl',\cl''}Gr(A_{\cl''})$$
where $f_{\cl,\cl',\cl''}\in\ca$ is $0$ for all but finitely many $[\cl'']\in B_G$.
We see that the family of elements $(Gr(A_\cl))_{\cl\in B_G}$ span a module that
is closed under multiplication.
Hence the class of simple perverse sheaves $(A_\cl)_{\cl\in B_G}$ on $G$ 
is of the kind described in the Introduction. Note that this class contains the character 
sheaves on $G$; moreover it contains only few non-character sheaves (compared to character 
sheaves).

\subhead 2.5. Example\endsubhead
In this subsection we assume that $G=PGL_2(\kk)$. 
There are exactly three strata of $G$: the set $Y_{rs}$ of regular semisimple elements; the 
set $Y_{ru}$ of regular unipotent elements; the set $\{1\}$.

We have $Y_{rs}=Y_{T,T}$ where $T$ is a maximal torus of $G$. Let $\cs^1(T)$ be the subcategory of $\cs(T)$ whose
objects are the $\ce\in\cs(T)$ which have rank $1$. For each $\ce\in\cs^1(T)$ we set $\cl_\ce=\p_!(\tce)$. If 
$\ce^{\ot 2}\not\cong\bbq$ then $\cl_\ce=\cl_{\ce^*}$ is irreducible (of rank $2$). If $\ce^{\ot2}\cong\bbq$ then 
$\cl_\ce\cong\cl'_\ce\op\cl''_\ce$ where $\cl'_\ce,\cl''_\ce$ are local systems of rank $1$ on $Y_{rs}$ such that
$\cl'_\ce$ extends to a local system on $G$ and $\cl''_\ce$ does not.
Now $B_{Y_{rs}}$ consists of $[\cl_\ce]=[\cl_{\ce^*}]$ (with $\ce\in\cs^1(T)$, $\ce^{\ot2}\not\cong\bbq$) and of 
$\cl'_\ce,\cl''_\ce$ (with $\ce\in\cs^1(T)$, $\ce^{\ot2}\cong\bbq$). If $\ce,\ce_1\in\cs^1(T)$ 
we have 

$[\cl_\ce][\cl_{\ce_1}]=[\cl_{\ce\ot\ce_1}]+[\cl_{\ce\ot\ce_1^*}]$. 
\nl
If in addition we have $\ce^{\ot2}\cong\bbq$ then 

$[\cl'_\ce][\cl_{\ce_1}]=[\cl''_\ce][\cl_{\ce_1}]=[\cl_{\ce\ot\ce_1}]$. 
\nl
If in addition we have $\ce^{\ot2}\cong\bbq$, $\ce_1^{\ot2}\cong\bbq$ then

$[\cl'_\ce][\cl'_{\ce_1}]=[\cl''_\ce][\cl''_{\ce_1}]=[\cl'_{\ce\ot\ce_1}]$,

$[\cl'_\ce][\cl''_{\ce_1}]=[\cl''_{\ce\ot\ce_1}]$. 
\nl
We have $Y_{ru}=Y_{G,S}$ where $S=Y_{ru}$ and $B_{Y_{ru}}$ consists of $[\bbq^{ru}]$ where 
$\bbq^{ru}$ is the local system $\bbq$ on $Y_{ru}$. We have 
$[\bbq^{ru}][\bbq^{ru}]=v[\bbq^{ru}]$.

We have $\{1\}=Y_{G,S}$ where $S=\{1\}$ and $B_{\{1\}}$ consists of $[\bbq^1]$ where $\bbq^1$ is the local system
$\bbq$ on $\{1\}$. We have $[\bbq^1][\bbq^1]=v^3[\bbq^1]$.

Let $\ce\in\cs^1(T)$, $\ce^{\ot2}\not\cong\bbq$. We write $A_\ce$ instead of $A_{\cl_\ce}$. 
Let $\ce\in\cs^1(T)$, $\ce^{\ot2}\cong\bbq$. We write $A'_\ce$ instead of $A_{\cl'_\ce}$ and 
$A''_\ce$ instead of $A_{\cl''_\ce}$. We write $A_{ru}$ instead of $A_{\bbq^{ru}}$ and $A_1$ 
instead of $A_{\bbq^1}$. We have

$Gr(A_\ce)=Gr(A_{\ce^*})=[\cl_\ce]+v\i[\bbq^{ru}]+(v^{-1}+v^{-3})[\bbq^1]$
if $\ce^{\ot2}\not\cong\bbq$,

$Gr(A'_\ce)=[\cl'_\ce]+v\i[\bbq^{ru}]+v^{-3}[\bbq^1]$,
$Gr(A''_\ce)=[\cl''_\ce]+v\i[\bbq^1]$ if $\ce^{\ot2}\cong\bbq$,

$Gr(A_{ru})=[\bbq^{ru}]+v^{-2}[\bbq^1]$, 

$Gr(A_1)=[\bbq^1]$.

From these formulas and from the multiplication table with respect to the basis $B_G$ we can easily compute the 
product of any two elements in the basis $\tB_G$ as an $\ca$-linear combination of elements in $\tB_G$. For 
example, if $\ce,\ce_1\in\cs^1(T)$, $\ce^{\ot2}\not\cong\bbq$, $\ce_1^{\ot2}\not\cong\bbq$, 
$(\ce\ot\ce_1)^{\ot2}\not\cong\bbq$, $(\ce\ot\ce^*_1)^{\ot2}\not\cong\bbq$, we have

$Gr(A_\ce)Gr(A_{\ce_1})=Gr(A_{\ce\ot\ce_1})+Gr(A_{\ce\ot\ce_1^*})-v\i Gr(A_{ru})+vGr(A_1)$.


\widestnumber\key{BBD}
\Refs
\ref\key{\BV}\by D.Barbasch and D.Vogan\paper Primitive ideals and orbital integrals in complex exceptional groups
\jour J.Algebra\vol80\yr1983\pages350-382\endref
\ref\key{\BFO}\by R.Bezrukavnikov, M.Finkelberg and V.Ostrik\paper Character D-modules via Drinfeld center
\finalinfo(in preparation)\endref
\ref\key{\BBD}\by A.Beilinson, J.Bernstein and P.Deligne\paper Faisceaux pervers\jour Ast\'erisque\vol100\yr1982
\endref
\ref\key{\GI}\by V.Ginzburg\paper Admissible modules on a symmetric space\jour Ast\'erisque\vol173-174\yr1989
\pages199-255\endref 
\ref\key{\KL}\by D.Kazhdan and G.Lusztig\paper Representations of Coxeter groups and Hecke algebras\jour Invent.
Math.\vol53\yr1979\pages165-184\endref
\ref\key{\KLL}\by D.Kazhdan and G.Lusztig\paper Schubert varieties and Poicar\'e duality\inbook Proc.Symp.Pure 
Math.\vol36\publ Amer.Math.Soc.\yr1980\pages185-203\endref
\ref\key{\MV}\by I.Mirkovi\'c and K.Vilonen\paper Characteristic varieties of character sheaves\jour Invent.Math.
\vol93\yr1988\pages405-418\endref 
\ref\key{\CS}\by G.Lusztig\paper Character sheaves,I\jour Adv. Math.\vol56\yr1985\pages193-237\moreref II\vol57\yr
1985\pages226-265\moreref III\vol57\yr1985\pages266-315\moreref IV\vol59\yr1986\pages1-63\moreref V\vol61\yr1986
\pages103-155\endref
\ref\key{\RA}\by G.Lusztig\paper Rationality properties of unipotent representations\jour J.Algebra\vol258\yr2002
\pages1-22\endref
\ref\key{\PAR}\by G.Lusztig\paper Parabolic character sheaves,I\jour Moscow Math.J.\vol4\yr2004\pages153-179
\endref
\ref\key{\CDG}\by G.Lusztig\paper Character sheaves on disconnected groups,I\jour Represent. Th. (electronic)\vol7
\yr2003\pages374-403\moreref II\vol8\yr2004\pages72-124\moreref III\vol8\yr2004\pages125-144\moreref IV\vol8\yr
2004\pages145-178\moreref Errata\vol8\yr2004\pages179-179\moreref V\vol8\yr2004\pages346-376\moreref VI\vol8\yr
2004\pages377-413\moreref VII\vol9\yr2005\pages209-266\moreref VIII\vol10\yr2006\pages314-352\moreref IX\vol10\yr
2006\pages353-379\moreref X\finalinfo(in preparation)\endref
\endRefs
\enddocument